\documentclass[12pt]{amsart}
\usepackage{amsmath,amscd,amssymb,latexsym}
\usepackage{graphicx}
\input{epsf}
\newtheorem{theorem}{Theorem}
\newtheorem{lemma}{Lemma}
\newtheorem{cor}{Corollary}
\newtheorem{prop}{Proposition}
\newtheorem{definition}{Definition}

\begin{document}

\title{The Yang-Mills Measure in the $SU(3)$ Skein Module}
\author{Charles Frohman}
\address{Department of Mathematics \\ The University of Iowa \\ Iowa-City, Iowa 52242}
\email{frohman@math.uiowa.edu}

\author{Jianyuan K. Zhong}
\address{Department of Mathematics and Statistics \\
California State University Sacramento\\ 6000 J Street \\Sacramento, CA 95826}
\email{kzhong@csus.edu }

\date{February 19, 2004}
\maketitle
\begin{abstract} Let $A\ne 0$ be a complex number that is not a root
  of unity.  Let $M$ be a compact smooth oriented $3$-manifold, 
the $SU(3)$-skein space of $M$, $S_A(M)$, is the vector space 
over $\mathbb{C}$ generated by framed oriented links (including 
framed oriented trivalent graphs in $M$) quotient by the 
$SU(3)$-skein relations due to Kuperberg. For a closed, orientable 
surface $F$, we construct a local diffeomorphism invariant trace on $S_A(F\times I)$.
\end{abstract}


\section{Introduction}


Throughout this paper, three manifolds and surfaces will be
compact and oriented. A {\bf framed oriented trivalent graph} is a
space that is homeomorphic to a closed regular neighborhood of an
oriented
trivalent graph embedded in an orientable surface, along with an
embedding of that oriented  graph in the space. As these are oriented, 
each edge of the
graph carries a direction. In diagrams, we will just draw the
graph and the reader can imagine its regular neighborhood running
parallel to the graph in the plane of the paper. We always have
the same ``side'' of the neighborhood facing up. By a {\bf framed
oriented link} in a three-manifold $M$ we mean an embedding of
such a space in $M$. Two framed oriented links are equivalent if
there is an isotopy of $M$ taking one to the other that preserves
the orientations of the edges. We will also
work with {\bf relative} framed oriented links. These are  graphs that also
have some monovalent 
vertices but
they are exactly the points of intersection of the graph with the
boundary of the manifold. Of course, the framed graph intersects
the boundary of $M$ in arcs so that each arc has a monovalent
vertex in its interior.

 Let $A\neq 0$ be a complex number so that if $A$ is a root of unity then
 $A=\pm 1$.
 When $A\neq \pm 1$, we define
\[ [n]=\frac{A^{3n}-A^{-3n}}{A^3-A^{-3}},\]
and when $A=\pm 1$, we define $[n]=n$. Let
$[n]!=[n][n-1]\cdots [1].$ Finally,  \[ \left[\begin{matrix}
n\\k\end{matrix}\right]=\frac{[n]!}{[k]![n-k]!}.\]

If $M$ is a compact oriented three-manifold let $\mathcal{L}$ be
the set of equivalence classes of framed oriented links so that
all the vertices are {\em sources} or {\em sinks}. That is, at
each vertex either all  three edges  point  in or they all point out. It is
worth noting that  the empty link is included in this collection.
Let $\mathbb{C}\mathcal{L}$ denote the vector space having
$\mathcal{L}$ as a basis. Let $R_A$ be the subspace of
$\mathbb{C}\mathcal{L}$ spanned by the following five skein relations from Kuperberg \cite{K}:

\begin{itemize}

\item {\em positive crossing}

\[\raisebox{-8pt}{\includegraphics{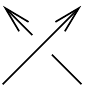}}- A^2 \raisebox{-8pt}{\includegraphics{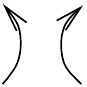}}
+A^{-1}\raisebox{-8pt}{\includegraphics{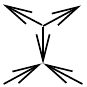}}\]

\item {\em negative crossing}

\[\raisebox{-8pt}{\includegraphics{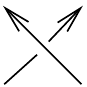}}- A^{-2} \raisebox{-8pt}{\includegraphics{cup.pdf}}
+A\raisebox{-8pt}{\includegraphics{tvfig.pdf}}\]

\item {\em square}

\[
\raisebox{-8pt}{\scalebox{.6}{\includegraphics{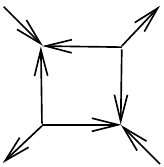}}}-\raisebox{-7pt}
{\includegraphics{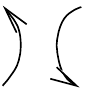}}\raisebox{1pt}{-}\raisebox{-7pt}{\includegraphics{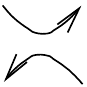}}\]

\item {\em bubble}

\[
\raisebox{-12pt}{\scalebox{.66}{\includegraphics{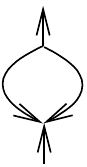}}}-[2]\raisebox{-8pt}{\includegraphics{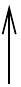}}\]

\item {\em trivial component}

\[ \raisebox{-6pt}{\includegraphics{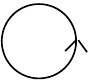}}D-[3]D\]
\end{itemize}

\begin{definition} The $SU(3)$ skein space of $M$ at $A$,  denoted by
$S_A(M)$, is the quotient space $\mathbb{C}\mathcal{L}/R_A$.
\end{definition}

There is another skein relation that can be easily derived from these,
which indicates change of framing.

\[ \raisebox{-16pt}{\scalebox{0.5}{\includegraphics{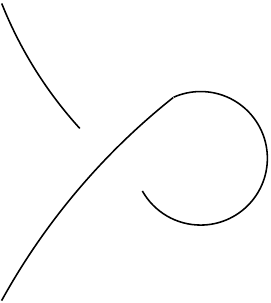}}}=A^8\quad
\raisebox{-10pt}{\includegraphics{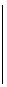}}\ .\]

It is convenient to note that the skein relations for the two
crossings, the change of framing and the trivial component generate
all the skein relations. In the case that $A=\pm1$, this skein module
has been studied by Adam Sikora \cite{S}. At these values of $A$,
crossings are irrelevant as the two crossing relations reduce to show
that the skeins are equal. Consequently there is a well defined
product structure. If $\gamma$ and $\beta$ are two framed oriented
graphs whose vertices are sources and sinks in $M$, we perturb them so
that they are disjoint and take their union as the product. The product structure on
$S_{\pm 1}(M)$ is the one induced by this. Sikora constructs a
natural homomorphism from $S_{\pm 1}(M)$ onto the $SL_3(\mathbb{C})$ characters
of $\pi_1(M)$ . The kernel of this homomorphism is the nilradical of $S_{\pm
1}(M)$.

We will also use relative skein spaces. For such, choose a
collection of  arcs in the boundary of $M$ along with a sign $[+]$ or $[-]$
for each arc.  The relative skein space is  the
vector space spanned by equivalence classes of relative framed
oriented links whose vertices are sources and sinks that intersect
the boundary of $M$ in those arcs so that if the sign of the arc is $[+]$ then
edge of the graph points into $M$ and if the sign of the is $[-]$
then the edge of the graph points out of $M$.

An especially important class of relative skein spaces are cylinders
over a disk, where we have indicated a family of arcs on the boundary
of the disk. You can think of assigning plusses and minuses to the
arcs, and either by enhancing the arguments of Sikora or imitating
the work of Kuperberg. The associated relative
module is isomorphic to $\mathrm{Inv}(V^p\otimes V^{*n})$, where $V$
is the fundamental representation of $SL_3\mathbb{C}$ and $V^*$ is its
dual, and $p$ is the number of positive arcs and $n$ the number of
negative arcs.

If $M=F\times [0,1]$, we  represent framed oriented links by
drawing oriented trivalent graphs with overcrossings and
undercrossings in $F$ and using the blackboard framing. In this
case $S_A(F\times [0,1])$ is an algebra. The multiplication is
defined by laying one skein over the other. To emphasize that the
algebra structure comes from the surface, we denote such a skein
space by $S_A(F)$.

In section 2, we recall some related results in $SU(3)$ skein and study the
$SU(3)$-skein modules of 
the solid torus $S^1\times D^2$, $S^1\times S^2$ and for the connected sum of
two $3$-manifolds. We list 
some results as follows. When $A\ne 0$ and $A$ is not a root of unity,

(1) $S_A(S^1\times D^2)$ has a countable basis indexed by the set of all ordered pairs of nonnegative integers.

(2) $S_A(S^1\times S^2)=\mathbb{C}\emptyset$, i.e., $S_A(S^1\times S^2)$ is generated by the empty framed link.

(3) $S_A(M_1\# M_2)\cong S_A(M_1)\otimes S_A(M_2)$. This says that the
$SU(3)$-skein module 
of the connected sum of two $3$-manifolds is isomorphic to the tensor product
of the $SU(3)$-skein 
modules of the manifolds.

(4) From (2) and (3), we conclude that the $SU(3)$-skein module of the
connected sum of 
$g$ copies of $S^1\times ^2$ is also generated by the empty skein.
$$S_A(\#^g S^1\times S^2)=\mathbb{C}\emptyset.$$

In sections 3, 4, we define and study the Yang-mills measure in a handlebody and on a closed surface. 

\section{Basics in $SU(3)$-skein theory}

\subsection{Related results from Ohtsuki and Yamada \cite{OY}---Magic Elements}
\begin{definition}
A magic element of type $(n,0)$ is inductively
defined by the following formula:
\[\raisebox{-4pt}{\includegraphics{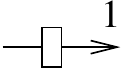}}=\raisebox{-3pt}{\includegraphics{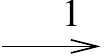}}\]
\[\raisebox{-4pt}{\includegraphics{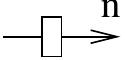}}=\raisebox{-13pt}{\includegraphics{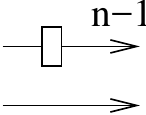}}-\frac{[n-1]}{[n]}\quad\raisebox{-14pt}{\includegraphics{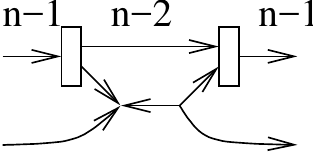}}
\]
\end{definition}
The following diagrams are called a left-Y and a right-Y:
\[\raisebox{-8pt}{\includegraphics{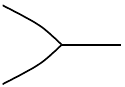}}, \quad
\raisebox{-8pt}{\includegraphics{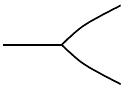}}\]

{\em Properties of the magic element of type $(n,0)$:  }

(1) When attached a left-Y to the right side or a right-Y to the
left side, the magic element of type $(n,0)$ vanishes.

(2) The magic element of type $(n,0)$ absorbs any magic elements
of type $(m,0)$ with $m\leq n$.
\[\raisebox{-10pt}{\includegraphics{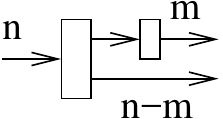}}=\raisebox{-10pt}
{\includegraphics{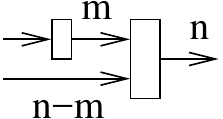}}=\raisebox{-4pt}{\includegraphics{nzero.pdf}}\]
\begin{definition}
 A magic element of type $(n,m)$ is defined by the
following formula:
\[\raisebox{-15pt}{\includegraphics{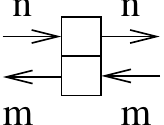}}=\sum_{k=0}^{\min(n,m)}(-1)^k
\frac{\left[\begin{matrix}n\\
k\end{matrix}\right]\left[\begin{matrix}m\\k\end{matrix}\right]}{\left[\begin{matrix}
n+m+1\\k\end{matrix}\right]}\quad\raisebox{-19pt}{\includegraphics{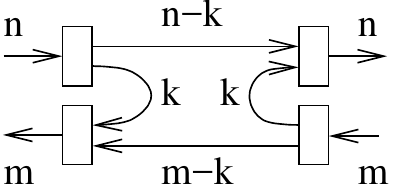}}\]
\end{definition}

We illustrate the left-U and right-U as follows:

\[\raisebox{-8pt}{\includegraphics{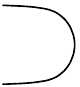}},\quad  \raisebox{-8pt}{\includegraphics{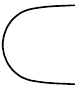}}\]
{\em Properties of the magic element of type $(n,m)$}:
When attached a left-Y or a left-U to the right side, or attached
a right-Y or a right-U to the left side, the magic element of type
$(n,m)$ vanishes.

\subsection{Coloring a trivalent graph with magic elements}

The coloring of an oriented edge $\includegraphics{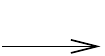}$ by
a pair of nonnegative integers $(n,m)$ is by replacing the edge in
the graph by the magic element of type $(n, m)$:
\[\raisebox{-8pt}{\includegraphics{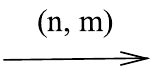}}=\raisebox{-16pt}{\includegraphics{nm.pdf}}\]
Then
\[\raisebox{-8pt}{\includegraphics{6edgenm.pdf}}=\raisebox{-8pt}{\includegraphics{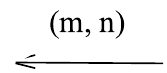}}\]

A vertex with acceptable labels becomes a triad with three edges
colored by (acceptable) nonnegative integer pairs $(n,m)$,
$(r,s)$ and (p,q). Here we illustrate a triad with indicated choice of orientations of the edges:
\[\raisebox{-12pt}{\includegraphics{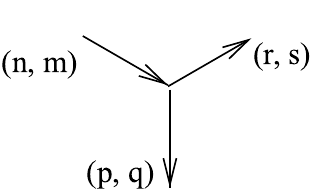}}\]
A triad represents a skein element in the relative skein space of
the disk with $(n+s+q)$ input points and $(m+r+p)$ output points.
Since there are possibly many different ways that strands can
intertwine in the middle, a triad with edges colored by $(n,m)$, $(r,s)$ and (p,q) is not uniquely
defined. Therefore we introduce a label by an $a_*$ inside a
circle to represent a specific intertwining of strands in the
middle of a triad and indicate it as
\[\raisebox{-12pt}{\includegraphics{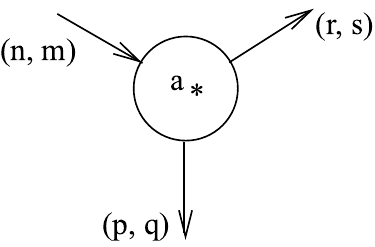}}\]

In the case that $A$ is either $\pm 1$ or not a root of unity, the
skein module can be interpreted in terms of invariant tensors in the
representation theory of $U_q(sl_3))$ (or $U(sl_3)$ when $A=\pm
1$). This can be seen in the works of Kuperberg \cite{K}, Kuperberg-Khovanov \cite{KK},
or Sikora \cite{S}. We summarize  some
conclusions of their work that we need for this development. Let $V$ be the
fundamental representation of $U_q(sl_3)$ and let $V^*$ be its
dual. Let $V_{p,q}$ be the highest weight irreducible representation
in $V^{\otimes p}\otimes V^{* \otimes q}$. There are invariant tensors
in $V \otimes V \otimes V$ and $V^* \otimes V^* \otimes V^*$ that
correspond to the two trivalent vertices. There are invariant
pairings
$V\otimes V^* \rightarrow \mathbb{C}$, and $V^*\otimes V \rightarrow
\mathbb{C}$, that can be used to ``stitch''  two trivalent vertices
together  along an edge so that each ``web'' in a disk, that is an embedded graph
with trivalent vertices in the interior of the disk and monovalent
vertices
on the boundary, so that the trivalent vertices are sources and sinks,
corresponds to an invariant tensor in the tensor product of copies of
$V$ and $V^*$ corresponding to choosing a basepoint on the boundary of
the disk and keeping track of the arrows going in and out as you go
around the disk. Modding out by the skein relations corresponding to
removing trivial simple closed curves, bubbles and four-sided regions
yields a vector space that is isomorphic to the space of
invariants. There is a further refinement where you group bunches of
edges together to form a clasped web space. Some of these ``clasped
web'' spaces
correspond to  $S_A(D^2, (n,m), (r,s), (p,q))$, the relative skein
space generated by the triads attached with the three magic
elements of types $(n,m)$, $(r,s)$ and $(p,q)$. We
call $(n,m), (r,s), (p,q)$ an admissible (acceptable) coloring of
a vertex if $S_A(D^2, (n,m), (r,s), (p,q))\ne 0$. From now on, we will only consider admissible
triads.  There are three
conclusions that we need to draw from this work.
\begin{enumerate}
\item The relative skein $S_A(D^2, (n,m), (r,s), (p,q))$ is up to
  cyclic permutation canonically isomorphic to
  $\mathrm{Inv}(V_{n,m}\otimes V_{r,s} \otimes V_{p,q})$.
Since $\mathrm{dim}(V_{m,n})\leq 9mn $, we see that there is a polynomial
$p(m,n,r,s,p,q)$ so that \[\mathrm{dim}(S_A(D^2, (n,m), (r,s), (p,q)))\leq p(m,n,r,s,p,q).\]
\item The pairing $S_A(D^2,(n,m),(r,s),(p,q))\otimes
  S_A(D^2,(q,p),(s,r),(m.n)) \rightarrow \mathbb{C}$ corresponding to
gluing the two disks together to form a sphere is nondegenerate as it
  corresponds to pairing  $\mathrm{Inv}(V_{m,n}\otimes V_{r,s}\otimes V_{p,q})$ with its
  dual. We can express this pairing in a diagram theoretic fashion as:
\[\raisebox{-12pt}{\includegraphics{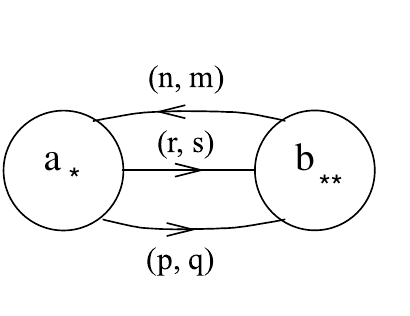}}\]
Note that $S_A(D^2)={\mathbb{C}}\emptyset$, so it induces a
pairing $< , >$ on \[S_A(D^2, (n,m), (r,s), (p,q)) \otimes
S_A(D^2, (q,p), (s,r), (m,n))\to \mathbb{C}\] by
\[<\alpha, \beta> =f(\alpha,\beta)\]
where $f(\alpha,\beta)$ is the complex number such that
$\alpha\otimes \beta =f(\alpha,\beta)\emptyset$ in $S_A(D^2)$. In
this sense, we have ${S_A(D^2, (n,m), (r,s), (p,q))}^* = S_A(D^2,
(q,p), (s,r), (m,n))$.

Following \cite{BB}, we choose bases $\{a_i\}$ for the space
$S_A(D^2, (n,m), (r,s), (p,q)$ and dual bases $\{b_i\}$ for
$S_A(D^2, (n,m), (r,s), (p,q))^{*}$,
so that
\[\raisebox{-16pt}{\includegraphics{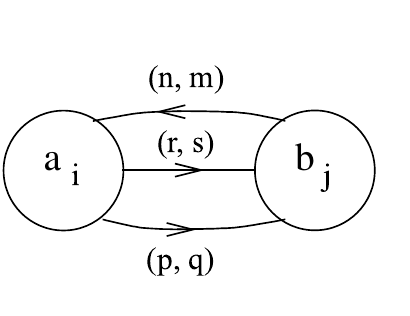}}\raisebox{24pt}{$=\delta_i^j$}\]

When one of the labels is $(0,0)$ then depending on the direction of the other
two arrows, the other two labels are either the same $(n,m)$ or its dual. In this case the skein
of the disk with two clasps is one dimensional so everything can be written as a multiple of 
the skein $s$ obtained by filling in the disk with straight lines from one clasp to the other.
Surprisingly, our choice of normalization leads to the peculiar realization that if our dual bases
are chosen as
$a=\alpha s$ and $b=\beta s$ then $\alpha*\beta=\frac{1}{\Delta_{m,n}}$.

\item The skein  space $S_A(D^2,(m,n),(r,s),(p,q),(u,v))$ is isomorphic
  to the result of ``stitching'' the sum below along the $(k,l)$ and $(l,k)$ factors.
\[\oplus S_A(D^2,(m,n),(r,s),(k,l))\otimes S_A(D^2,(l,k),(p,q),(u,v)),\]
where the sum is over all $(k,l)$ so that 
  $\mathrm{Inv}(V_{m,n}\otimes V_{r,s} \otimes V_{k,l})$ is nonzero
  and $\mathrm{Inv}(V_{l,k}\otimes V_{r,s}\otimes V_{u,v})$ is nonzero.
Using the dual bases chosen above we get the fusion formula \cite{BB}:

\[\raisebox{-24pt}{\includegraphics{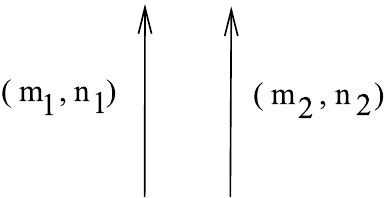}}=\sum \Delta_{k,l}\raisebox{-86pt}{\includegraphics{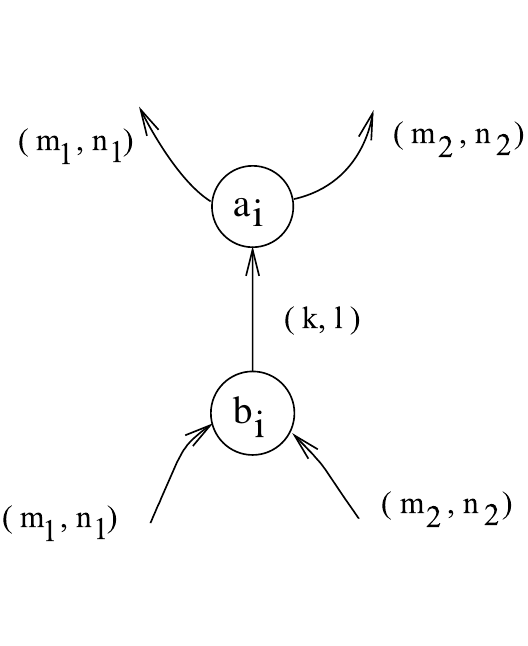}}.\]where the sum is over all admissible triples $(m_1,n_1), (m_2,n_2), (k,l)$ and dual bases $a_i$ and $ b_i$.

\end{enumerate}

Let $P_{n,m}$ be the closure of the magic element of type $(n,m)$
in the solid torus $S^1\times D^2$:
\[P_{n,m}=\raisebox{-44pt}{\scalebox{0.6}{\includegraphics{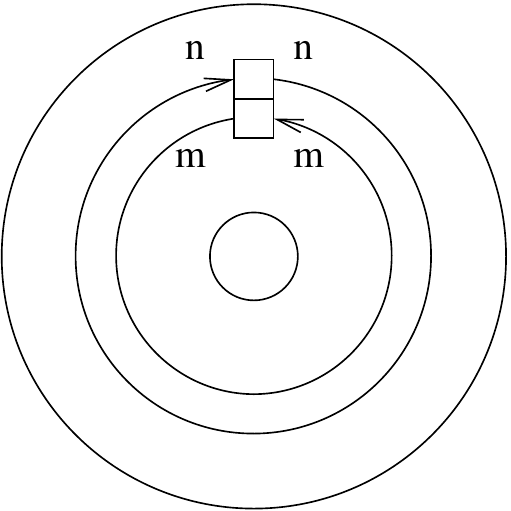}}}\]

As we show the solid torus, it is the cylinder over $S^1\times [0,1]$. The
inclusion
of $S^1 \times [0,1]$ into the plane $\mathbb{R}^2$ induces a corresponding
inclusion
of cylinders.
Let $\Delta_{n,m}$ be the complex scalar multiple of $P_{n,m}$ by
writing $P_{n,m}=\Delta_{n,m}\emptyset$ in
$S_A({\mathbb{R}}^2)=\mathbb{C}\emptyset$ by including $S_A(S^1\times
[0,1]$ into $S_A({\mathbb{R}}^2)$.

The following identities hold:

(1) From \cite{OY}, when $m,n$ are nonnegative integers, $\Delta_{n,m}=[n+1][m+1][n+m+2]/[2]$. When at least one of $m,n$ is a negative integer, we define $\Delta_{n,m}=0$. Note that if $A$ is a real number and is not a root of unity, then $\Delta_{n,m}\geq 1$ for all nonnegative integers $m,n$.

(2)
\[\raisebox{-86pt}{\includegraphics{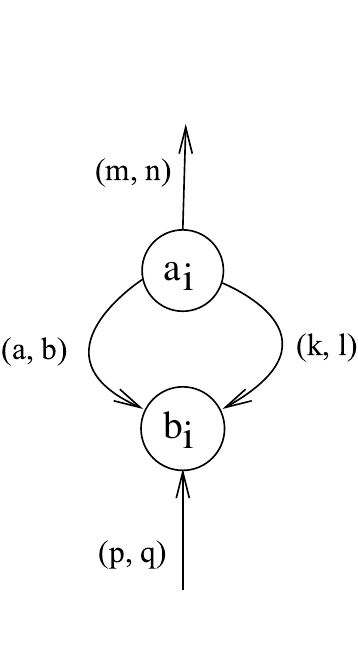}}=\frac{\delta_{m,n}^{p,q}}{\Delta_{m,n}}\quad\raisebox{-20pt}{\includegraphics{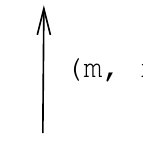}}\]
where $a_i$ and $b_i$ are dual bases for the admissible triple
$(m,n),\ (a,b), \ (k,l)$.

\begin{proof} (I) When $(m,n)\ne (p,q)$, we can prove that the skein element on the left hand side is zero. This eventually follows from the non-convexity of the basis of clasped web space of Kuperberg \cite{K}. To explain, we consider the clasped web space $W(C)$ where $C$ is given by the sequence $[(+\cdots +)_n(-\cdots -)_m(+\cdots +)_p(-\cdots -)_q]$ (the subscript indicates the number of plusses and minuses in the sequence. When $(m,n)\ne (p,q)$, the clasped web space $W(C)$ is zero, as there will be a minimal cut path with lower weight separating the clasps $[(+\cdots +)_n(-\cdots -)_m]$ and $(+\cdots +)_p(-\cdots -)_q]$ which causes convex clasps.

(II) When $(m,n)= (p,q)$, we prove the skein element on the left hand side is a scalar multiple of the magic element of type $(m,n)$. 

According to the Lemma 3.3 of \cite{OY}, if $D$ is a diagram in the disk with $2n+2m$ boundary points with neither biangles nor squares as shown:
\[D=\raisebox{-10pt}{\scalebox{0.6}{\includegraphics{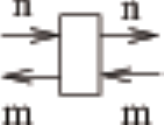}}}\]
Then $D$ satisfies at least one of the following three conditions:

(i) There is a left-Y or a left-U attached to the left side,

(ii) There is a right-Y or a right-U attached to the right side,

(iii) The diagram $D$ is $n+m$ parallel lines from the left side to the right side.

We can rewrite the middle part of the skein diagram of the given identity as a linear sum of diagrams fitting in the above lemma, then all three possible cases will give a multiple of the magic element of type $(m,n)$. Note that cases (i) and (ii) will contribute $0$ when attached to the magic element of type $(m,n)$.

(III) To find the scalar multiple, we close both sides, only when $a_i$ and $b_i$ are dual bases, the closure of the left hand side is nonzero and equals $1$, the closure on the right hand side of the magic element of type $(m,n)$ contributes $\Delta_{m,n}$. Therefore the scalar multiple is $\frac{1}{\Delta_{m,n}}$, the identity holds.
\end{proof}

(3) Let (k,l) be a pir of nonnegative integers such that the triples $(m_3,n_3),\ (m_1, n_1), \ (k,l)$ and $(m_2,n_2),\ (m_4,
n_4), \ (k,l)$ are admissible, then the collection of
elements
\[\includegraphics{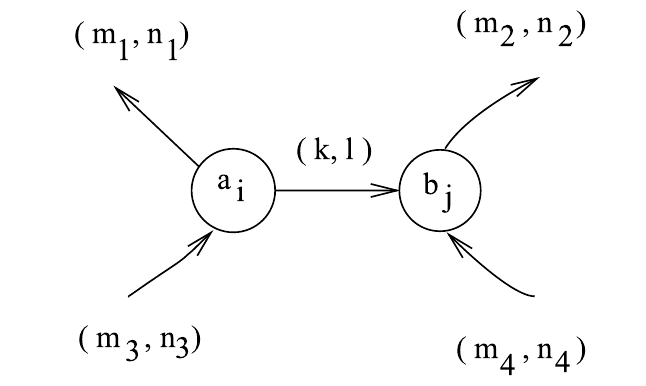}\]
 over all such $(k,l)$ and all basis elements $a_i$, $b_j$ of the corresponding triads, forms a basis for the skein
space $S_A(D^2, (m_3,n_3),\ (m_1, n_1), \ (m_2,n_2),\ (m_4,
n_4))$.

On the other hand, over all pairs of nonnegative integers $(g,h)$ such that  the  triples $(m_1,n_1),\ (m_2, n_2), \ (g,h)$ and $(m_3,n_3),\ (m_4,
n_4), \ (g,h)$ are admissible and basis elements $c_p$, $d_q$ of the corresponding triads, the collection of elements
\[\includegraphics{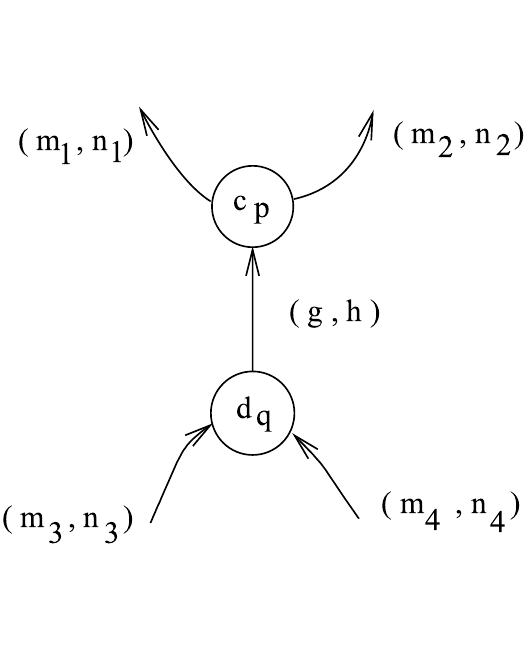}\]
also forms a basis for the skein space $S_A(D^2, (m_3,n_3),\ (m_1, n_1), \ (m_2,n_2),\ (m_4, n_4))$.

(4) There is a change of bases on $S_A(D^2, (m_3,n_3),\ (m_1, n_1), \ (m_2,n_2),\ (m_4, n_4))$:
\[\raisebox{-54pt}{\includegraphics{mn1234.pdf}}\]
\[=\sum \left
\{\begin{matrix} (m_3,n_3),\ (m_1, n_1)\ (g,h)\\   (m_2, n_2)\
(m_4, n_4)\ (k,l)\end{matrix} \ a_i, b_j, c_p,
d_q\right\}\raisebox{-80pt}{\includegraphics{mm1234.pdf}}\] where
the summation is over all bases
$c_p, d_q$ and admissible triples $(m_1,n_1),
(m_2,n_2), (g,h)$ and $(m_3,n_3), (m_4,n_4), (g,h)$.

Similarly, we have a pairing on $S_A(D^2, (m_1,n_1),\ (m_2, n_2),
\ (m_3,n_3),\ (m_4, n_4))\otimes S_A(D^2, (n_1,m_1),\ (n_2, m_2),
\ (n_3,m_3),\ (n_4, m_4))\to {\mathbb{C}}$ induced by the bilinear
form of attaching skein elements along the boundary of $D^2$
through the inclusion into $S_A(D^2)$.
\[\includegraphics{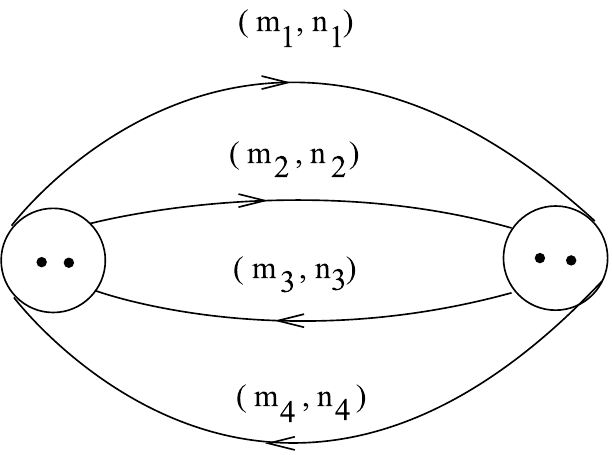}\]

\begin{definition}
We define \[\mathrm{Tet}\left\{ \begin{matrix} ((m_3,n_3),\ (m_1, n_1)\ (g,h)\\ (m_2, n_2)\ (m_4, n_4)\ (k,l)\end{matrix}\ a_i, b_j,
c_p, d_q\right\}\] to be the complex multiple of writing the
following skein element as a multiple of the empty skein
$\emptyset$ in $S_A(D^2)=\mathbb{C}\emptyset$.
\[\includegraphics{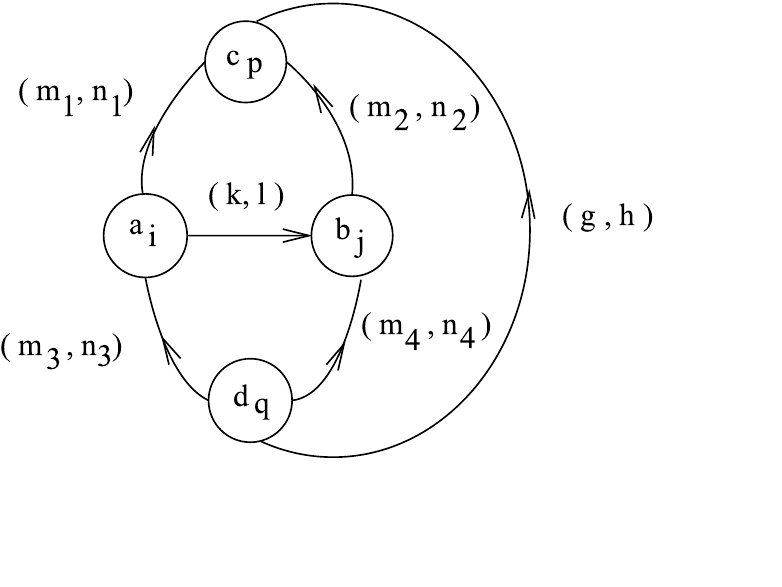}\]
\end{definition}

\begin{theorem}
\[\left\{ \begin{matrix} (m_3,n_3),\ (m_1, n_1)\ (g,h)\\   (m_2, n_2)\ (m_4, n_4)\ (k,l)\end{matrix}\ a_i,\ b_j,\ c_p,\ d_q\right\}
=\mathrm{Tet}\left\{ \begin{matrix} (m_3,n_3),\ (m_1, n_1)\ (g,h)\\
(m_2, n_2)\ (m_4, n_4)\ (k,l)\end{matrix}\ a_i,\ b_j,\ c_p,\
d_q\right\}\Delta_{g,h}\]
\end{theorem}

The proof of this theorem is similar to the computation in \cite{KL}.

\begin{theorem} If  $A>0$ and $A$ is not a root of unity,
\[\Biggl|\ \mathrm{Tet}\left\{ \begin{matrix} (m_3,n_3),\ (m_1, n_1)\ (g,h)\\
(m_2, n_2)\ (m_4, n_4)\ (k,l)\end{matrix}\ a_i,\ b_j,\ c_p,\
d_q\right\}\ \Biggr|\leq
\frac{1}{\sqrt{\Delta_{k,l}\Delta_{g,h}}}\]
\end{theorem}
\begin{proof}
The key is using the change of bases identity twice.
\[\raisebox{-54pt}{\includegraphics{mn1234.pdf}}\]
\[=\sum \left
\{\begin{matrix} (m_3,n_3),\ (m_1, n_1)\ (g,h)\\   (m_2, n_2)\
(m_4, n_4)\ (k,l)\end{matrix} \ a_i, b_j, c_p,
d_q\right\}\raisebox{-80pt}{\includegraphics{mm1234.pdf}}\]
\[=\sum \alpha_* \raisebox{-54pt}{\includegraphics{mn1234.pdf}}\]
where
\[\alpha_*= \left
\{\begin{matrix} (m_3,n_3),\ (m_1, n_1)\ (g,h)\\   (m_2, n_2)\
(m_4, n_4)\ (k,l)\end{matrix} \ a_i, b_j, c_p, d_q\right\} \left
\{\begin{matrix} (m_4, n_4)\ (m_3, n_3)\ (k,l)\\   (m_1, n_1)\
(m_2, n_2)\ (g,h)\end{matrix} \ a_i, b_j, c_p,
d_q\right\} \] 
Therefore,
\[\sum \left
\{\begin{matrix} (m_3,n_3),\ (m_1, n_1)\ (g,h)\\   (m_2, n_2)\
(m_4, n_4)\ (k,l)\end{matrix} \ a_i, b_j, c_p, d_q\right\} \left
\{\begin{matrix} (m_3,n_3),\ (m_1, n_1)\ (g,h)\\   (m_2, n_2)\
(m_4, n_4)\ (k,l)\end{matrix} \ a_i, b_j, c_p, d_q\right\}=1\]
i.e.,
\[\sum \mathrm{Tet}^2\left\{ \begin{matrix} (m_3,n_3),\ (m_1, n_1)\ (g,h)\\
(m_2, n_2)\ (m_4, n_4)\ (k,l)\end{matrix}\ a_i\ b_j\ c_p\
d_q\right\}\Delta_{k,l}\Delta_{g,h}=1\] We observe that each term
in the summation is positive, hence \[0\leq \mathrm{Tet}^2\left\{ \begin{matrix} (m_3,n_3),\ (m_1, n_1)\ (g,h)\\
(m_2, n_2)\ (m_4, n_4)\ (k,l)\end{matrix}\ a_i\ b_j\ c_p\
d_q\right\}\Delta_{k,l}\Delta_{g,h}\leq 1\] The result follows.
\end{proof}

\begin{cor}
When $A>0$ and $A$ is not a root of unity,
\[\Biggl|\ \mathrm{Tet}\left\{ \begin{matrix} (m_3,n_3),\ (m_1, n_1)\ (g,h)\\
(m_2, n_2)\ (m_4, n_4)\ (k,l)\end{matrix}\ a_i,\ b_j,\ c_p,\
d_q\right\}\ \Biggr|\leq
\frac{1}{\sqrt{\Delta_{k,l}}}\]
and
\[\Biggl|\ \mathrm{Tet}\left\{ \begin{matrix} (m_3,n_3),\ (m_1, n_1)\ (g,h)\\
(m_2, n_2)\ (m_4, n_4)\ (k,l)\end{matrix}\ a_i,\ b_j,\ c_p,\
d_q\right\}\ \Biggr|\leq
\frac{1}{\sqrt{\Delta_{g,h}}}\]
\end{cor}
\begin{proof}
This follows from the fact that $\Delta_{g,h}\geq 1$ and $\Delta_{k,l}\geq 1$ for all nonnegative integer pairs $(g,h)$ and $(k,l)$. 
\end{proof}
\subsection{Some fundamental examples}

Assume $A$ is not a root of unity.
\begin{prop} \cite{OY}
$S_A(S^1\times D^2)$ has a basis given by the collection
$\{P_{n,m}\ |\ n,m $ $ \textrm{are\ nonnegative\ integers}\}$.
\end{prop}

\begin{theorem}\label{two}
\[S_A(S^1\times S^2)=\mathbb{C}\emptyset.\]
\end{theorem}
\begin{proof}
As $S^1\times S^2$ can be obtained from the solid torus $S^1\times
D^2$ by adding a $2$-handle. There is an epimorphism
$S_A(S^1\times D^2)\to S_A(S^1\times S^2)$ induced by embedding
$S^1\times D^2$ into $S^1\times S^2$. Adding a $2$ handle results
in adding relations to the generators. Here we prove it suffices
to consider only the following sliding relation:
\[\raisebox{-20pt}{\includegraphics{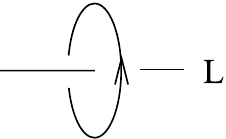}}=\raisebox{-20pt}{\includegraphics{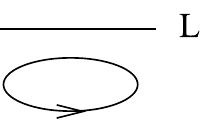}}\]
the equality holds in $S_A(S^1\times S^2)$, where $L$ is any skein
element in $S_A(S^1\times D^2)$. We only need to consider the
sliding relation on the generators $P_{n,m}$ of $S_A(S^1\times
D^2)$. From \cite{OY}, we have the following skein relation:
\[\raisebox{-20pt}{\includegraphics{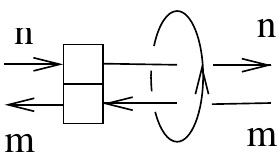}}=(A^{4n+2m+6}+A^{-2n+2m}+A^{-2n-4m-6})\raisebox{-20pt}{\includegraphics{nm.pdf}}\]
While the product of $P_{n,m}$ with a trivial component is
$[3]P_{n,m}=(A^6+A^{-6}+1)P_{n,m}$, we conclude that
$(A^{4n+2m+6}+A^{-2n+2m}+A^{-2n-4m-6}-A^6-A^{-6}-1)P_{n,m}=0$ 
in $S_A(S^1\times S^2)$. On the other hand, it's easy to obtain the following skein relation similar to the above:
\[\raisebox{-20pt}{\includegraphics{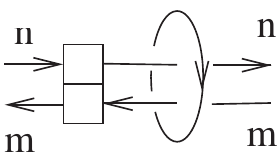}}=(A^{4m+2n+6}+A^{-2m+2n}+A^{-2m-4n-6})\raisebox{-20pt}{\includegraphics{nm.pdf}}\]
Then $(A^{4m+2n+6}+A^{-2m+2n}+A^{-2m-4n-6}-A^6-A^{-6}-1)P_{n,m}=0$ in
$S_A(S^1\times S^2)$.
 When $A$ is not a root of unity and $m,n$ are not both zero, one can prove that $(A^{4n+2m+6}+A^{-2n+2m}+A^{-2n-4m-6}-A^6-A^{-6}-1)$ and $(A^{4m+2n+6}+A^{-2m+2n}+A^{-2m-4n-6}-A^6-A^{-6}-1)$ are not both zero.
Therefore $P_{n,m}=0$ in $S_A(S^1\times S^2)$ when $n, m$ are not both
zero. While $P_{0,0}=\emptyset$ doesn't involve the sliding relation,
it survives. Hence $S_A(S^1\times S^2)=<\emptyset>$.
\end{proof}

\begin{theorem}
\[S_A(M_1\# M_2)=S_A(M_1)\otimes S_A(M_2)\]
\end{theorem}
The proof follows the same outline as \cite{GZ1}.

\begin{cor}
\[S_A(\#^g S^1\times S^2)=\mathbb{C}\emptyset.\]
\end{cor}

\section{The Yang-Mills measure in $S_A(F\times [0,1])$ with $\partial
F\ne \emptyset$}

Let $F$ be a compact oriented surface and $I=[0,1]$. Assume $A$ is not a root of unity. We denote the
skein algebra of $S_A(F\times I)$ by $S_A(F)$ to emphasize that
the algebra structure depends on $F$. Notice that $F\times I$ is a
handlebody. If you choose a family $K$ of proper arcs on $F$ that
cut it down to a disk, then $K\times I$ is a family of disks that
cut $F\times I$ into a ball. The double of $F\times I$, denoted by
$D(F\times I)$, is the result of gluing two copies of $F\times I$
together using the identity map on their boundary. The disks
$K\times I$ in each copy are glued together to form a system of
spheres in $D(F\times I)$ that cut it down to a punctured ball.
Therefore $D(F\times I)$ is homeomorphic to a connected sum of
copies of $S^1\times S^2$. From the Preliminaries in 2.4,
$S_A(D(F\times I))={\mathbb{C}}\emptyset$. This induces a linear
functional ${\mathcal{YM}}: S_A(F\times I)\to {\mathbb{C}}$ by the
inclusion of $F\times I$ into $D(F\times I)$, i.e., if $\alpha \in
S_A(F\times I)\subset S_A(D(F\times I))$, we can write
$\alpha=f(\alpha)\emptyset$ for some complex number $f(\alpha)$ in
$S_A(D(F\times I))$, then ${\mathcal{YM}}(\alpha)=f(\alpha)$ .

\begin{prop}
\[{\mathcal{YM}}(\alpha * \beta)={\mathcal{YM}}(\beta *\alpha)\]
\end{prop}

\begin{proof}
If you remove $\partial F \times I$ from $F\times I$, then the
double of the resulting object is homeomorphic to ${\mathring{
F}}\times S^1$ where the product structure coincides with the
product $\mathring{ F}\times I$ on each half. Also
\[F\times I\subset \mathring{
F}\times S^1 \subset D(F\times I),\] we can perturb
representations of $\alpha$ and $\beta$ so that they miss
$\partial F$. Now it is clear that $\alpha *\beta$ and $\beta *
\alpha$ are isotopic in $\mathring{ F}\times S^1$.
\end{proof}

Let $\Gamma$ be an oriented trivalent spine of $F$. An admissible
coloring of $\Gamma$ is by attaching the magic elements of types
$(m_i, n_i)$ along the edges and acceptable labels at the
vertices, such a coloring of $\Gamma$ is an element of $S_A(F)$.

\begin{theorem}
The admissible colorings of $\Gamma$ form a spanning set for $S_A(F)$.
\end{theorem}
\begin{proof} Let ${\mathcal{H}}_g$ be the handlebody $F\times I$, let $D$ be a separating meridian disk of ${\mathcal{H}}_g$:
\[\raisebox{-3mm}{\includegraphics{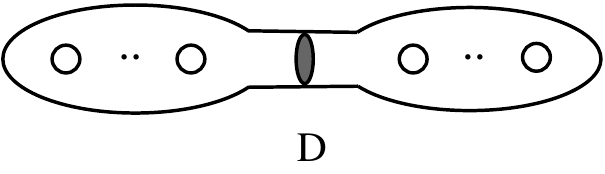}}\]

Let $V_D=[-1, 1]\times D$ be a regular neighborhood of $D$ in
${\mathcal{H}}_g$, $V_D$ can be projected into a disk $D_p=[-1,
1]\times [0, 1]$.

Let $\alpha$ be a framed link in ${\mathcal{H}}_g$ in general position to $V_D$, let $\alpha'=\alpha\cap
V_D$,
\[\alpha'=\raisebox{-10mm}{\includegraphics{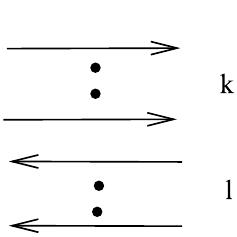}}\]
In the next lemma, we show that we can write $\alpha'$ as a linear
sum of skein elements in $V_D$ which have the magic elements in
the middle such as the following:
\[\raisebox{-3mm}{\includegraphics{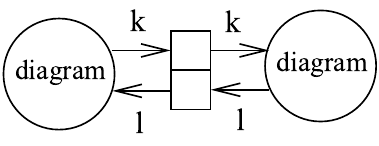}}\]
where $k,l$ are some nonnegative integers.

By induction on the number of separating meridian discs of
${\mathcal{H}}_g$, we can write $\alpha$ as a linear sum of skein
elements in $S_A({\mathcal{H}}_g)$ which have the magic elements
at the regular neighborhood of each separating disk. Note that
such an element corresponds to an admissible coloring of the
trivalent spine of ${\mathcal{H}}_g$. Therefore, the admissible
colorings of $\Gamma$ spans $S_A({\mathcal{H}}_g)$.
\end{proof}
\begin{lemma}\label{one}
An element in the relative skein module of the cylinder with $n$
parallel strands going to the right and $m$ parallel strands going
to the left can be written as a linear sum of elements in the form
\[\raisebox{-3mm}{\includegraphics{6nm.pdf}}\]
with $(k+l)\leq (m+n)$.
\end{lemma}
\begin{proof}
We proceed by induction on $n+m$.

(1) When $n+m=1$, i.e, $n=1, m=0$ or $m=1,n=0$, it is trivial.

When $n+m=2$, there are three cases: (i) $n=m=1$; (ii) $n=2, l=0$;
(iii) $m=2, n=0$. We illustrate cases (i) and (ii) by the
following. Case (iii) is similar to case (ii).
\[\raisebox{-8mm}{\includegraphics{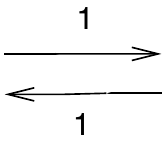}}=\raisebox{-8mm}{\includegraphics{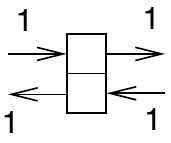}}+\frac{[1]}{[2]}\raisebox{-5mm}{\includegraphics{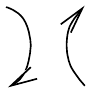}}\]
\[\raisebox{-8mm}{\includegraphics{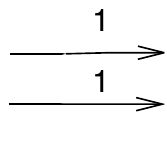}}=\raisebox{-6mm}{\includegraphics{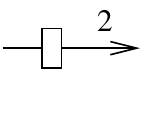}}+\frac{[1]}{[2]}\raisebox{-4mm}{\includegraphics{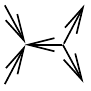}}\]

(2) Assume the result is true for $(n+m)\leq k$ for some natural
number $k$, we prove the result is true for the case
$(n+m)=(k+1)$. By induction, the result is true for all
nonnegative integers $n$, $m$.

When $n+m=k+1$, without loss of generality, we can assume that
$m\geq 1$. By the induction assumption, the result is true for
$n+(m-1)=k$, i.e., we can consider the part with $n$ parallel
strands going to the right and $(m-1)$ parallel strand going to
the left and write it as a linear sum of elements of the assumed
form with the size of the magic elements in the middle of $\leq
k$; now it suffices to prove that the following element is a
linear sum of the assumed form,
\[\raisebox{-3mm}{\includegraphics{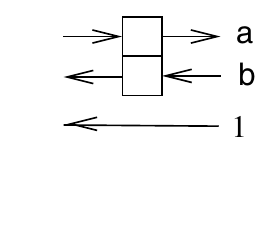}}\]
where $(a+b)= k$.

 By the definition of the magic element of type $(a,b)$,
\[\raisebox{-38pt}{\includegraphics{6nm1.pdf}}=\sum_{i=0}^{\min(a,b)}(-1)^i
\frac{\left[\begin{matrix}a\\
i\end{matrix}\right]\left[\begin{matrix}b\\i\end{matrix}\right]}{\left[\begin{matrix}
a+b+1\\i\end{matrix}\right]}\quad\raisebox{-40pt}{\includegraphics{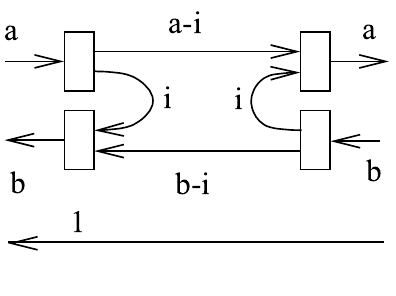}}\]
Note that all terms in the summation corresponding to $i\geq 1$
intersect the separating disk at no more than $k$ times, so by the
induction assumption, these can be written as the sum of the
assumed forms. Therefore we only need to worry about the first
term which corresponds to $i=0$, i.e.,

\[\raisebox{-12pt}{\includegraphics{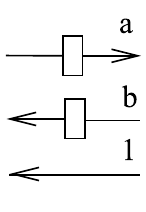}}\]
Now use the definition of the magic element of type $(0, b)$,
\[\raisebox{-24pt}{\includegraphics{6mn.pdf}}=\raisebox{-30pt}{\includegraphics{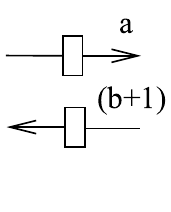}}+\frac{[b]}{[b+1]}\raisebox{-30pt}{\includegraphics{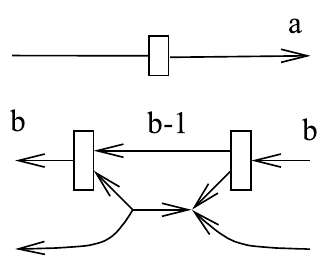}}\]
Similarly, the second term intersects the separating disk at
$(a+b)=k$ times, so we only need to consider the first term on the
right hand side, again we use the definition of the magic element
of type $(a, b+1)$:
\[\raisebox{-35pt}{\includegraphics{6m1n.pdf}}=\raisebox{-45pt}{\includegraphics{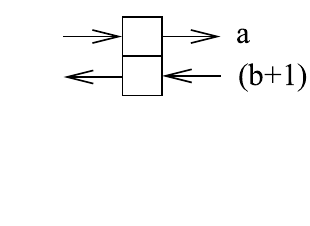}}-\sum_{i=1}^{\min(a,b)}(-1)^i
\frac{\left[\begin{matrix}a\\
i\end{matrix}\right]\left[\begin{matrix}b+1\\i\end{matrix}\right]}{\left[\begin{matrix}
a+(b+1)+1\\i\end{matrix}\right]}\raisebox{-55pt}{\includegraphics{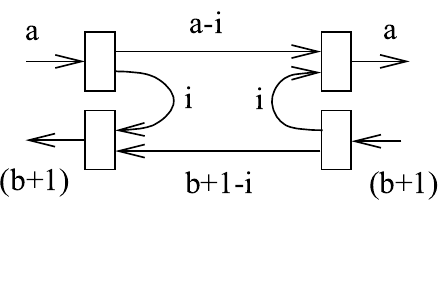}}\]
We observe again that each term in the summation has $i\geq 1$, so
the skein elements intersect the separating disk at $\leq k$
times; by the induction assumption, these can be written as a
linear sum of elements which have the magic elements at the middle
of size $\leq k$. Therefore we conclude that the result is true
for $n+m=k+1$.
\end{proof}
\begin{lemma}\label{three}
If $k$ is a properly embedded arc in $F$, then $S_A(F)$ is spanned by
trivalent colored graphs so that each graph intersects $k$ in one
transverse point of intersection in the interior of one of its edges.
\end{lemma}

The preceeding lemmas give a method of computing the Yang-Mills measure. Let 
$k$ be a proper arc in $F$. By Lemma \ref{one} any skein can be written as a 
colored trivalent graph intersecting $k$ in a single point of transverse 
intersection in the interior of an edge. By the argument in Theorem
\ref{two} the Yang-Mills measure of such a graph is zero unless that edge 
carries the label $(0,0)$. So, choose a system of proper arcs that cut $F$ 
into a family of disks. Using Lemma \ref{three} repeatedly write the skein 
as a colored trivalent graph that intersects each of the arcs at most once 
in a point of transverse intersection in the interior of an edge. Throw out 
all terms where a graph carries a nonzero label on an edge intersecting one 
of the arcs. Next erase the edges and renormalize to take into account the 
peculiarity of the normalization.  Finally, evaluate the invariant of the 
remaining skeins in the disks.

We formalize this:

\begin{prop} {\bf Locality} \label{locality}
Let $F$ be a compact oriented surface and let $k$ be a proper arc. Let $F'$ be 
the result of cutting $F$ along $k$. If $s$ is a skein that is represented by 
a sum of colored trivalent graphs that each intersects $k$ in at most a single 
point of transverse intersection in the interior of an edge. Let $s_0$ be sum 
of all terms where the  graph is either disjoint from $k$ or intersects in an edge labeled $(0,0)$. The skein $s_0$ corresponds to a skein $s'_0$ in the image of the inclusion of $S_A(F') \rightarrow S_A(F)$. Then,
\[ \mathcal{YM}(s)=\mathcal{YM}(s'_0).\] \end{prop}

\qed

\section{The Yang-Mills measure on a closed surface}
In this section $F$ will be a closed oriented surface of genus greater than
$1$. Further we suppose that $A$ is a positive real number not equal to $1$.
We prove that there is a linear functional
\[ \mathcal{YM}:S_A(F) \rightarrow \mathbb{C}.\]
that is a trace in the sense that
$\mathcal{YM}(\alpha*\beta)=\mathcal{YM}(\beta*\alpha)$ for all $\alpha, \beta
\in S_A(F)$.

If we remove an open disk with nice boundary from $F$ we get a compact surface
with one boundary component $F'$ which is a subsurface of $F$. The inclusion
map $i: F' \rightarrow F$ induces a surjective map,
\[ S_A(F') \rightarrow S_A(F).\] 
Let $\partial_{(m,n)}$ be the skein that is the result of coloring
a framed knot that is parallel to $\partial F$ with the $(m,n)$ type magic element.
We define $\mathcal{YM}:S_A(F) \rightarrow \mathbb{C}$. If $\alpha \in S_A(F)$
then choose $\alpha'$ to be a skein in $S_A(F')$ that gets mapped onto
$\alpha$ by the inclusion. We denote the Yang-Mills measure on $F'$ by
$\mathcal{YM}_F'$ and define,
\begin{equation} \mathcal{YM}(\alpha)=\lim_{N \rightarrow \infty} 
\sum_{(m,n)}^{m+n \leq N} \Delta_{m,n}
\mathcal{YM}_{F'}(\partial_{(m,n)}*\alpha').\tag{*}\end{equation} 

We need to prove two things. The first is that the series we gave above is
convergent and the second is to prove that is independent of the choice of
$\alpha'$.  To prove convergence we just need to prove it converges on a
collection of skeins that spans $S_A(F)$. Luckily we have such a family,
admissibly colored spines.  

Let $s_c$ denote a trivalent spine of a compact oriented surface $F$ with one
boundary component that has been colored admissibly. If
the spine has $v_i$ vertices and $e_j$ edges, then 
$$v_i=-2\chi(F);\quad e_j=-3\chi(F).$$ where $\chi(F)$ is the Euler characteristic of $F$. 

We need a global estimate on
$\mathcal{YM}(\partial_{m,n}*s_c)$. Let $(p_j,q_j)$ be the label on the $j$th
edge of $s_c$ and let $a_i$ be the skein in the $i$th vertex.  To compute this
we fuse along the handles, that look like this

\begin{picture}(130,112)
\hspace{1in} \scalebox{2}{\includegraphics{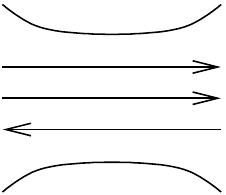}}
\put(0,72){$(m,n)$} \put(0,54){$(p_j,q_j)$} \put(0,36){$(m,n)$} 
\end{picture}

We compute by fusing one $(m,n)$ strand with the $(p_i,q_i)$ strand,
then fusing with the other $(m,n)$ strand, and throwing out everything that
the central edge is not labeled $(0,0)$. The number of terms is equal to the
multiplicity of $V_{(n,m)}$ in $V_{(m,n)}\otimes V_{(p_i,q_i)}$. We call the
set of pairs of skeins that appear in the two vertices along this edge,
$\mathcal{A}_j$.  We then need
to erase the $(0,0)$ edges, which entails dividing by $\Delta_{m,n}$ for each edge.
We are then just computing the value of sum of the products of some
tetrahedral coefficients.
The result is
\[ \sum_{\mathcal{A}_j}\prod_{v_i} \mathrm{Tet}
\left\{ \begin{matrix} (m, n) & (m,n) & (m,n) \\ (p_{i_1},q_{i_1})
    &(p_{i_2},q_{i_2}) &
    (p_{i_3},q_{i_3})  
\end{matrix}\ a_j, b_j, c_j,d_j\right\}\]
The $a_j, b_j,c_j,d_j$ are skeins in vertices coming from the fusions along the
edges. The number of terms is less than or equal to a polynomial evaluated on
the labels, and the size of each term in the sum is less than $\prod_i
\frac{1}{\sqrt{\Delta_{m,n}}}$. Therefore we have the following Proposition.

\begin{prop} There is a polynomial in variables $(m,n)$ that only depends on
  the colors assigned to the edges, $p(m,n)$ so that 
\[ \mathcal{YM}(\partial_{(m,n)}*s_c)\leq
\frac{p(m,n)}{\Delta_{m,n}^{-\chi(F)} },\]
where $\chi(F)$ is the Euler characteristic of $F$.\end{prop}

\begin{prop} The formula given by the equation (*) for the Yang-Mills measure converges.\end{prop} 
\proof The proof is by comparison with the series $\sum_{m,n}
\frac{p(m,n)}{\Delta_{m,n}^{1/2}}$. We know from its formula that 
$\Delta_{m,n}$ grows exponentially in $m$ and $n$. Hence the series we just
mentioned converges. By the estimate given by the previous proposition and since
the euler characteristic of $F'$ is $\leq -3$ we see that the terms
in the series for $\mathcal{YM}(s_c)$ are bounded in absolute value by the
series we just gave. \qed

The final step of the argument is to show that $\mathcal{YM}(\alpha)$ is independent of the choice of $\alpha'$. Since we can pass from any skein $\alpha'$ to any other skein $\alpha''$ that is sent to $\alpha$ under
$$S_A(F')\to S_A(F)$$
by handleslides.
By fusing, we can reduce this to check this is true for the result of sliding one string of a trivalent colored spine of $F'$ across the boundary disk. Without loss of generality, let $\alpha'=s_c$ be the trivalent spine of a compact oriented surface $F$ with one
boundary component that has been colored admissibly,
let $\partial_{m,n}$ be the framed knot corresponding to $\partial F$ oriented with
the boundary orientation from $F$ colored with the $(m,n)$ magic element.  

Locally $\partial_{(m,n)}*\alpha'$  looks like
\[\raisebox{-8pt}{\includegraphics{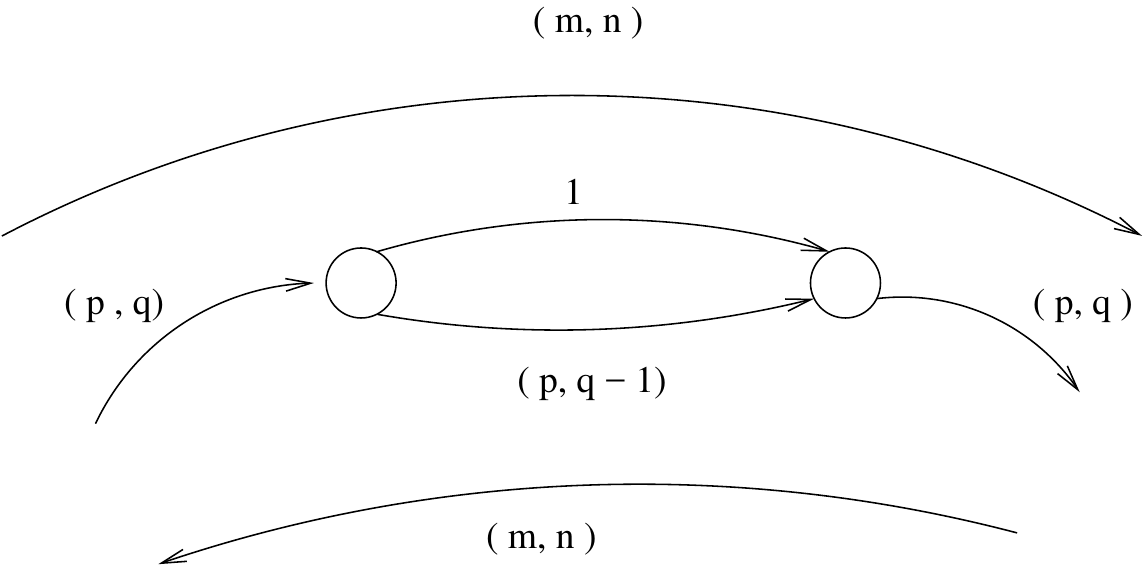}}\]
Let $\alpha''$ be the skein obtained from $\alpha'=s_c$ by sliding one strand over the added disk, locally the diagram $\partial_{(m,n)}*\alpha''$ looks like
\[\raisebox{-8pt}{\includegraphics{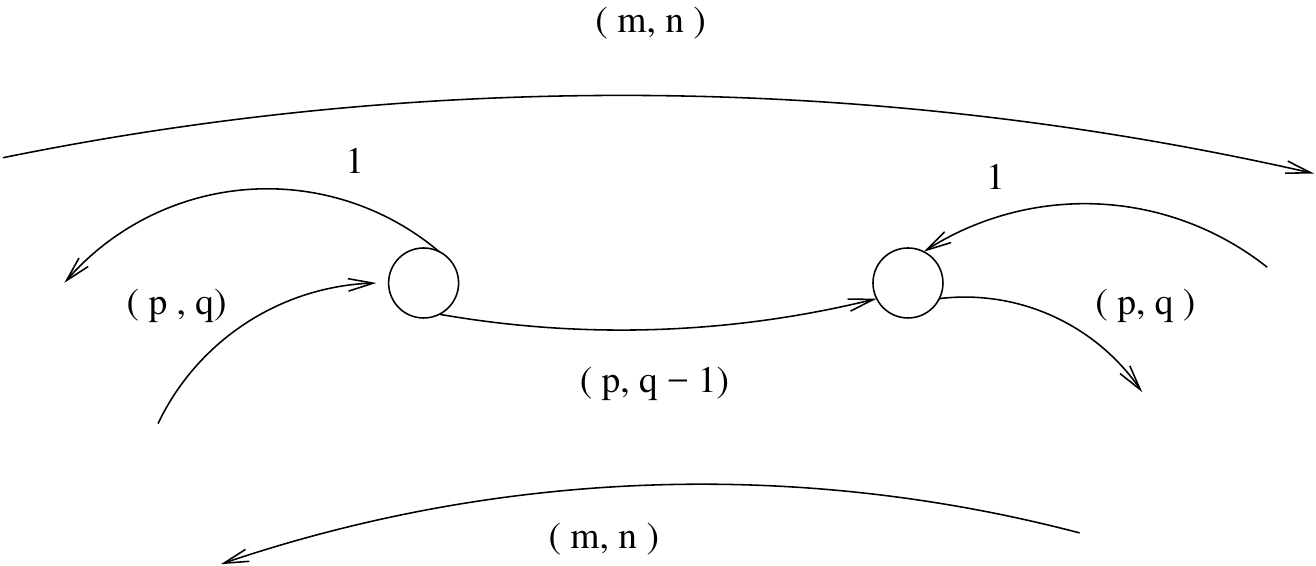}}\]

In the following, we will show that the Yang-Mills measure defined by the equation (*) is well-defined by proving
$$\lim_{N\to \infty}\sum_{(m,n)}^{m+n \leq N} \Delta_{m,n}
(\mathcal{YM}_{F'}(\partial_{(m,n)}*\alpha')-
\mathcal{YM}_{F'}(\partial_{(m,n)}*\alpha''))=0. $$

\begin{lemma}
$$\sum_{(m,n)}^{m+n \leq N} \Delta_{m,n}
((\partial_{(m,n)}*\alpha')-
(\partial_{(m,n)}*\alpha'')) = \sum_{(m,n)}^{m+n = N } \Delta_{m,n}(\Delta_{m+1, n}s_1-\Delta_{n, m+1}s_2)$$
where $s_1$ illustrated below is the skein element which is almost the same as $\partial_{(m,n)}*\alpha'$ except where it's shown, 
\[s_1=\quad \raisebox{-65pt}{\includegraphics{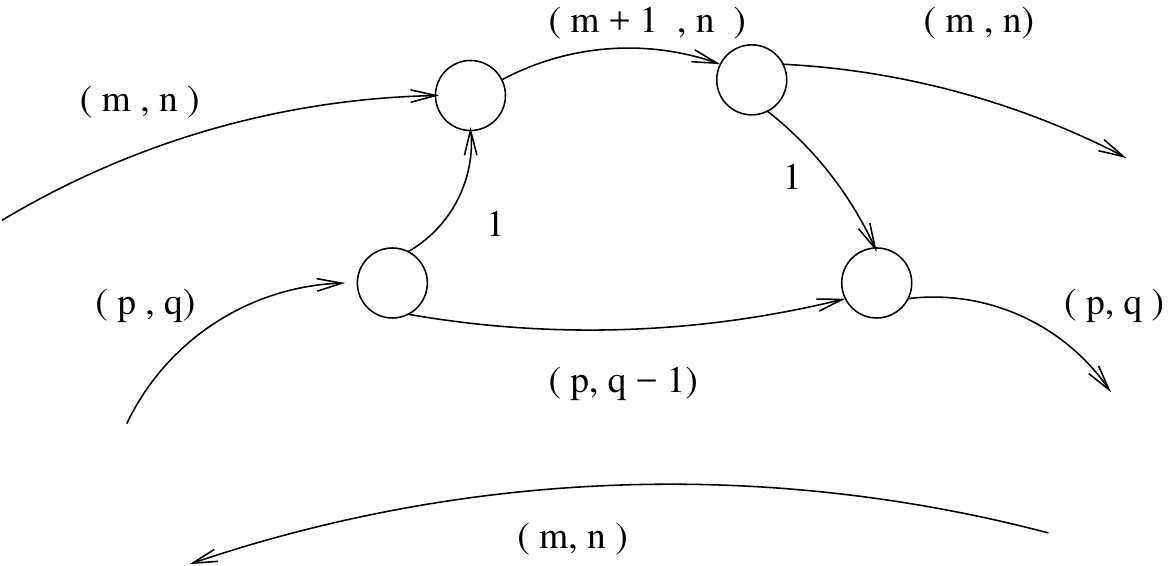}}\]
$s_2$ is the skein element which is almost the same as $\partial_{(m,n+1)}*\alpha''$ except where it's shown,
\[s_2=\quad \raisebox{-65pt}{\includegraphics{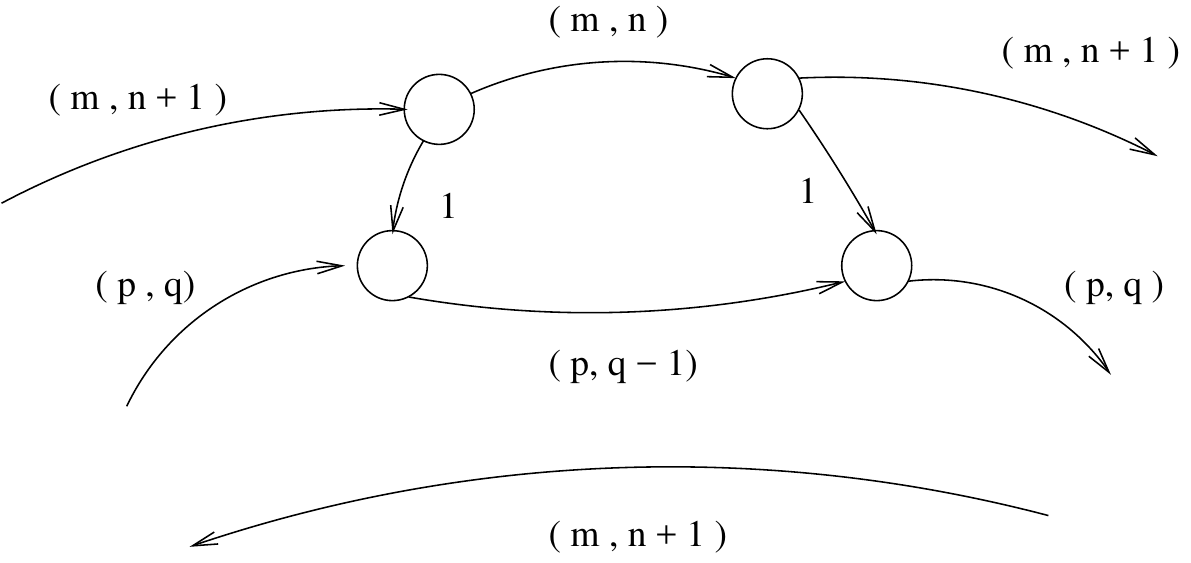}}.\]
\end{lemma}
The proof follows from the following lemma.

\begin{lemma}
$$\sum_{(m,n)}^{m+n \leq N} \Delta_{m,n}
(\quad \raisebox{-65pt}{\includegraphics{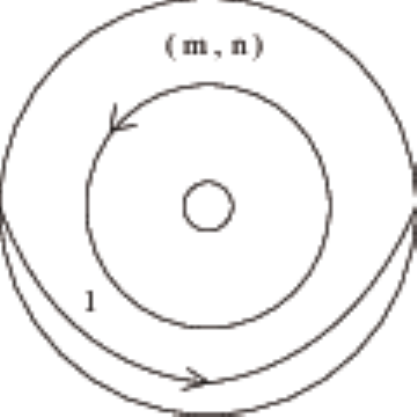}}-\quad \raisebox{-65pt}{\includegraphics{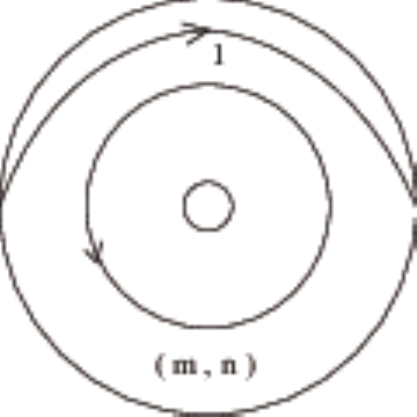}}\quad) = $$
$$\sum_{(m,n)}^{m+n = N } \Delta_{m,n}(\Delta_{m+1, n}\quad \raisebox{-65pt}{\includegraphics{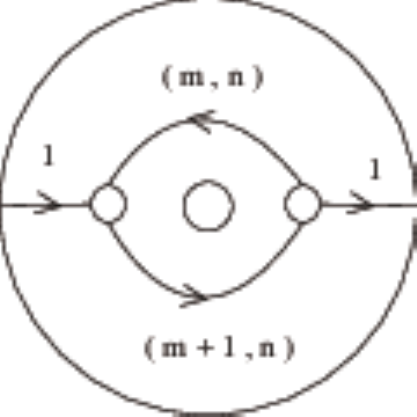}}-\Delta_{m, n+1}\quad \raisebox{-65pt}{\includegraphics{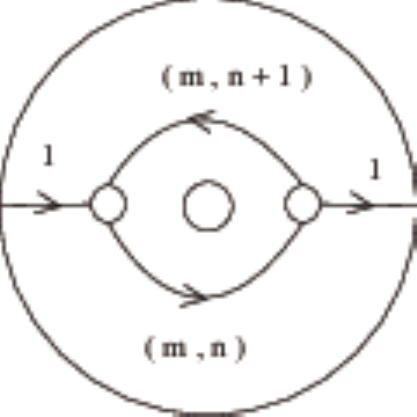}}\ )$$
\end{lemma}
\begin{proof}
We consider the tensor product of the magic element of type $(m,n)$ and the magic element of type $(1,0)$. From the representation theory, we have $$V_{m,n}\otimes V_{1,0}=V_{m,n-1}\oplus V_{m-1,n+1}\oplus V_{m+1,n}.$$ Similarly,
$$V_{m,n}\otimes V_{0,1}=V_{m-1,n}\oplus V_{m+1,n-1}\oplus V_{m,n+1}.$$
These give the corresponding fusion identities in the $SU(3)$-skein. When we apply these fusion identities to the left hand side of the identity in the lemma, almost all terms are canceled except the terms left on the right hand side.  
\end{proof}

\begin{theorem}
$$\lim_{N \rightarrow \infty} 
\sum_{(m,n)}^{m+n = N} \Delta_{m,n}\Delta_{m+1,n}
\mathcal{YM}_{F'}(s_1)=0$$
\end{theorem}
\begin{proof}
To compute the Yang-Mills measure of the skein $s_1$, we fuse to isolate the vertices. Notice that fusing $s_1$ will require two more cross cuts than that of $\partial_{(m,n)}*\alpha'$. After throwing out everything that the central edge is not labeled $(0,0)$ and erasing the $(0,0)$ edges, the Yang-Mills measure of $s_1$ is the product of 
\begin{equation}\sum_{\mathcal{A}'_j}\frac{1}{\Delta_{m,n}\Delta_{m+1,n}} \mathrm{Tet}
\left\{ \begin{matrix} (m, n) & (m,n) & (m+1,n) \\ (1,0)
    &(p,q-1) &
    (p,q)  
\end{matrix}\ a_j, b_j, c_j,d_j\right\}\mathrm \times \tag{I}\end{equation}
$${Tet}
\left\{ \begin{matrix} (m, n) & (m+1,n) & (m,n) \\ (1,0)
    &(p,q) &
    (p,q-1)    
\end{matrix}\ a'_j, b'_j, c'_j,d'_j\right\}$$
(the $a_i, b_j,c_j,d_j$ and $a'_i, b'_j, c'_j,d'_j$ are skeins in vertices coming from the fusions along the
edges.)

\noindent with the standard product of fusion on $\partial_{(m,n)}*\alpha'$,
\begin{equation}  \sum_{\mathcal{A}_j}\prod_{v_i\ \mathrm{of}\ \alpha'} \mathrm{Tet}
\left\{ \begin{matrix} (m, n) & (m,n) & (m,n) \\ (p_{i_1},q_{i_1})
    &(p_{i_2},q_{i_2}) &
    (p_{i_3},q_{i_3})  
\end{matrix}\ a_i, b_j, c_j,d_j\right\}\tag{II}\end{equation} 
First the product in (I) is less than or equal to 
$$\sum_{\mathcal{A}_j}\frac{1}{\Delta_{m,n}\Delta_{m+1,n}}\frac{1}{\sqrt{\Delta_{m,n}}\sqrt{\Delta_{m+1,n}}}\leq p'(m,n) \Delta_{m,n}^{-\frac{3}{2}}\Delta_{m+1,n}^{-\frac{3}{2}}$$
where $p'(m,n)$ is another polynomial depending on the colors $(m,n)$.

Secondly, by a previous proposition, there exist polynomials $p(m,n)$ in variables $(m,n)$ that depend only on the colors assigned to the edges so that the above standard product is less than or equal to
$$\frac{p(m,n)}{\Delta_{m,n}^{-\chi(F)}}.$$

Therefore 
$$\mathcal{YM}_{F'}(s_1)\leq p'(m,n)p(m,n) \Delta_{m,n}^{-\frac{3}{2}}\Delta_{m+1,n}^{-\frac{3}{2}}{\Delta_{m,n}^{\chi(F)}},$$
and
$$\lim_{N \rightarrow \infty} 
\sum_{(m,n)}^{m+n = N} \Delta_{m,n}\Delta_{m+1,n}
\mathcal{YM}_{F'}(s_1)=\lim_{N \rightarrow \infty} 
\sum_{(m,n)}^{m+n = N} p'(m,n)p(m,n) \Delta_{m,n}^{-\frac{1}{2}+\chi(F)}\Delta_{m+1,n}^{-\frac{1}{2}}=0$$
\end{proof}

\begin{cor}
$$\lim_{N \rightarrow \infty} 
\sum_{(m,n)}^{m+n = N} \Delta_{m,n}\Delta_{n,m+1}
\mathcal{YM}_{F'}(s_2)=0$$
\end{cor}

\begin{cor}
$$\lim_{N \rightarrow \infty}\sum_{(m,n)}^{m+n \leq N} \Delta_{m,n}
(\mathcal{YM}_{F'}(\partial_{(m,n)}*\alpha')-
\mathcal{YM}_{F'}(\partial_{(m,n)}*\alpha''))=$$
$$\lim_{N \rightarrow \infty} 
\sum_{(m,n)}^{m+n = N} \Delta_{m,n}(\Delta_{m+1,n}
\mathcal{YM}_{F'}(s_1)-\Delta_{n,m+1}
\mathcal{YM}_{F'}(s_2)=0$$
\end{cor}



\begin{thebibliography}{10}


\bibitem{AM98}
A.~K.~Aiston and H.~R.~Morton, {\em Idempotents of Hecke algebras
of type A}, J. of Knot Theory and Ram. 7 No 4 (1998), 463-487.

\bibitem{BB}
A.~Beliakova and C.~Blanchet, {\em Modular categories of types B,
C and D}, Comment. Math. Helv. 76 (2001) 467-500.


\bibitem{GZ1}
P.~Gilmer and J.~K.~Zhong, {\em The Homflypt skein module of a
connected sum of $3$-manifolds}, Algebraic and Geometric Topology,
Volume 1 (2001), 627-686.

\bibitem{KK}
M.~Khovanov and G.~Kuperberg, {\em Web bases for sl(3) are not dual canonical}, Pacific J. Math. 188 (1999), 129--153.

\bibitem{K}
G.~Kuperberg, {\em Spiders for rank $2$ Lie algebra}, Comm. Math. Phys. 180(1):109-151, 1996.

\bibitem{KL} L. H. Kauffman and S. Lins, {\em Temperley-Lieb
recoupling theory and invariants of $3$-manifolds}, Ann.\ of Math.\
Studies {\bf 143}, Princeton University Press, 1994.



\bibitem{OY}
T.~Ohtsuki and S.~Yamada, {\em Quantum $SU(3)$ Invariants of
$3$-manifolds Via linear Skein Theory}, Journal of Knot Theory and
Its Ramfications, Vol. 6, No. 3 (1997) 373-404.

\bibitem{S} A.~Sikora {\em Quantum $SU(n)$ Skein Theory} 

\bibitem{TW}
V.~Turaev and H.~Wenzl, {\em Quantum invariants of $3$-manifolds
associated with classical simple Lie algebras }, Internat. J.
Math., 4 (1993) 323-358.




\end{thebibliography}
\end{document}